\documentclass[11pt]{article}
\usepackage[T1]{fontenc}
\usepackage{lmodern,amsmath,amssymb,amsthm,microtype}
\usepackage[margin=1in]{geometry}
\usepackage[hidelinks]{hyperref}
\newtheorem{theorem}{Theorem}[section]

\newtheorem{corollary}[theorem]{Corollary}
\theoremstyle{definition}

\DeclareMathOperator{\Sing}{Sing}
\DeclareMathOperator{\Reg}{Reg}
\DeclareMathOperator{\vol}{vol}
\title{Ergodicity of geodesic flows on\\Aravinda--Farrell twisted doubles}
\author{Alberto Verjovsky}
\date{}
\begin{document}
\maketitle
\begin{abstract}
The announced rank-one ergodicity theorem of Rodriguez Hertz, Rodriguez Hertz, and Wu includes the nonpositively curved metrics considered here on Aravinda--Farrell twisted doubles. For these construction metrics, a direct geometric argument identifies the singular set. Their singular geodesics are precisely those tangent to the torus sections of finitely many flat product collars, so their initial vectors have Liouville measure zero. A compact simply connected modification region cannot contain a complete geodesic; this also treats the exotic connected-sum examples. Ergodicity then follows from the classical theorem on the regular set.
\end{abstract}



\section{The geometric criterion}

The work of Federico Rodriguez Hertz, Jana Rodriguez Hertz, and Weisheng Wu \cite{RHW} led me to examine the Aravinda--Farrell examples. Their announced theorem includes these examples. For the particular metrics considered here, the flat collars allow a short geometric verification of ergodicity and a description of the singular geodesics.

The metric constructions are due to Aravinda and Farrell \cite{AF09}, using Heintze's cusp flattening as recalled in \cite[Proposition 2.1]{AF09}, the warped-product formulas of Bishop and O'Neill \cite{BO}, and the local exotic-metric construction of Farrell and Jones \cite{FJ}. The dynamical argument belongs to the theory initiated by Hopf \cite{Hopf} and developed by Anosov, Sinai, and Pesin \cite{Anosov,AS,Pesin}. Rank is used in the sense of Ballmann, Brin, and Eberlein \cite{BBE}. 

Let $(N,g)$ be a closed connected smooth Riemannian manifold of dimension $d\geq2$. Write $SN$ for its unit tangent bundle, $\pi:SN\to N$ for the bundle projection, $g^t$ for its geodesic flow, and $m$ for normalized Liouville measure. Let $\gamma_v$ be the unit-speed geodesic with initial velocity $v$. If $d\sigma_p$ denotes probability spherical measure on $S_pN$, then
\[
 dm(v)=\frac{d\vol_g(p)}{\vol_g(N)}\,d\sigma_p(v),\qquad p=\pi(v).
\]
Ergodicity means that every measurable set invariant under $g^t$ has $m$-measure zero or one.

Assume that $g$ has nonpositive sectional curvature. A vector $v\in SN$ is \emph{singular} if the geodesic $\gamma_v$ admits a nonzero parallel normal Jacobi field. With $R$ denoting the Riemann curvature tensor, there is a field $J$ along $\gamma_v$ satisfying
\begin{equation}\label{jacobi}
 \langle J,\dot\gamma_v\rangle=0,\qquad
 \nabla_{\dot\gamma_v}J=0,\qquad
 R(J,\dot\gamma_v)\dot\gamma_v=0.
\end{equation}
Write $\Sing$ for the singular set and $\Reg=SN\setminus\Sing$. Vectors in $\Reg$ are called rank-one vectors; $N$ has rank one if $\Reg\neq\varnothing$.

The classical ergodicity theorem on the regular set states that for a closed connected rank-one manifold of nonpositive sectional curvature, the restriction of Liouville measure to $\Reg$ is ergodic. In particular,
\begin{equation}\label{criterion}
 m(\Sing)=0\quad\Longrightarrow\quad (g^t,m)\text{ is ergodic}.
\end{equation}
This is the consequence of the Hopf--Pesin theory needed here; see \cite{Pesin} and the explicit use of this implication in \cite[\S5, proof of Theorem 5.1(3)]{CLMT}.

\begin{theorem}\label{main}
Suppose that $N$ has nonpositive sectional curvature and contains finitely many pairwise disjoint compact product collars
\[
 C_i\cong T_i^{d-1}\times[a_i,b_i],\qquad
 g|_{C_i}=h_i+dr^2,\qquad a_i<b_i,
\]
where each $(T_i,h_i)$ is a flat torus. Suppose also that there are finitely many compact connected simply connected smooth domains $B_j$, pairwise disjoint and disjoint from all the $C_i$, such that every sectional curvature is strictly negative on the nonempty set
\[
 U=N\setminus\left(\bigcup_i C_i\cup\bigcup_j B_j\right).
\]
Either family may be empty. Then
\begin{equation}\label{singular}
 \Sing=\bigcup_i\{v\in SN|_{C_i}:\langle v,\partial_r\rangle=0\}.
\end{equation}
In particular, $m(\Sing)=0$ and the geodesic flow is ergodic.
\end{theorem}

\begin{proof}
If $v$ is singular, the parallel field in \eqref{jacobi} has constant positive length. At any time when $\gamma_v(t)\in U$, the Jacobi equation would imply
\[
 0=\langle R(J,\dot\gamma_v)\dot\gamma_v,J\rangle
   =K(\dot\gamma_v\wedge J)|J|^2<0.
\]
Thus $\gamma_v(\mathbb R)$ is contained in the union of the collars and the domains. These are finitely many disjoint compact sets. The connected image of $\mathbb R$ must therefore lie in one of them.

A complete geodesic cannot remain in a domain $B_j$. Indeed, the universal covering map $q:\widetilde N\to N$ restricts to a trivial covering over the simply connected domain $B_j$. Every component of $q^{-1}(B_j)$ is consequently compact. A geodesic contained in $B_j$ would lift to a complete geodesic contained in one such component. By the Cartan--Hadamard theorem, $\widetilde N$ is a complete simply connected manifold of nonpositive curvature, and each of its unit-speed geodesics satisfies
\[
 d_{\widetilde N}(\widetilde\gamma(s),\widetilde\gamma(t))=|s-t|.
\]
This is incompatible with containment in a compact set.

It follows that a singular geodesic remains in a collar $C_i$. In its product coordinates, write $\gamma_v(t)=(x(t),r(t))$. The geodesic equation gives $r''(t)=0$, hence
\[
 r(t)=r(0)+ct.
\]
Since $r(t)\in[a_i,b_i]$ for all $t\in\mathbb R$, we have $c=0$. Conversely, a geodesic tangent to a torus section remains in that section, and $\partial_r$ is a nonzero parallel normal Jacobi field along it. These statements hold also at the boundary sections, by smoothness of the ambient metric. This proves \eqref{singular}.

For every point in the interior of a collar, the vectors in \eqref{singular} form an equatorial $S^{d-2}$ in $S_pN\cong S^{d-1}$. They have spherical measure zero. The collar boundaries have base volume zero. The disintegration formula for $m$ and the finiteness of the collars give $m(\Sing)=0$. Finally, a vector based in $U$ has rank one, and \eqref{criterion} proves ergodicity.
\end{proof}

\section{Flattening the cusps}

Let $M$ be a compact connected smooth manifold with torus boundary such that $\operatorname{int}M$ admits a complete finite-volume real hyperbolic metric. An affine boundary diffeomorphism $f$ determines the twisted double
\[
 D_fM=M\cup_f M.
\]
Here affine refers to the affine structures on the cusp sections. The construction in \cite[Proposition 2.1 and Lemma 2.2]{AF09} changes the metric in the cusps to obtain flat product boundary collars, with prescribed boundary metrics up to sufficiently small scaling.

The flattening can be chosen with strictly negative sectional curvature everywhere before the product region. Choose a smooth function $q:\mathbb R\to[0,\infty)$, supported on $[0,2]$, positive on $(0,2)$, symmetric about $1$, and satisfying $\int_0^2q(u)\,du=1$. Set
\[
 \phi(s)=-s+\int_0^s(s-u)q(u)\,du.
\]
Then $\phi(s)=-s$ for $s\leq0$, while symmetry gives $\int_0^2u q(u)\,du=1$ and hence $\phi(s)=-1$ for $s\geq2$. Moreover,
\[
 \phi''\geq0,\qquad \phi'(s)<0\quad(s<2).
\]
For a flat metric $h$ on a cusp section, use
\[
 g=dt^2+\psi(t)^2h,\qquad
 \psi(t)=\exp\bigl(\phi(t-t_0)-t_0\bigr).
\]
This agrees with $dt^2+e^{-2t}h$ for $t\leq t_0$ and is a product for $t\geq t_0+2$. The warped-product curvature formulas \cite{BO} give
\begin{equation}\label{curvature}
 K_{\mathrm{tan}}=-\left(\frac{\psi'}{\psi}\right)^2=-(\phi')^2,
 \qquad
 K_{\mathrm{rad}}=-\frac{\psi''}{\psi}
              =-\bigl(\phi''+(\phi')^2\bigr).
\end{equation}
Derivatives of $\phi$ in \eqref{curvature} are evaluated at $t-t_0$. The curvature operator splits into the tangential two-forms and the forms containing $\partial_t$, with these two eigenvalues. Thus all sectional curvatures are negative for $t<t_0+2$ and vanish in the product region.

The other interpolation in \cite[Lemma 2.2]{AF09} connects two parallel flat metrics on the torus. Its metric is
\[
 A_\alpha=e^{-2\alpha t}G_t+\alpha^2dt^2,\qquad 1\leq t\leq2,
\]
where $G_t$ is a smooth family of parallel flat metrics, constant near the endpoints. For sufficiently large $\alpha$, all its sectional curvatures are strictly negative, as established in that lemma. These interpolations leave the hyperbolic core unchanged.

\begin{corollary}\label{affine}
Under the hypotheses of \cite[Theorem 1.2]{AF09}, the metric on $D_fM$ can be chosen so that its geodesic flow is ergodic for Liouville measure.
\end{corollary}

\begin{proof}
Choose the flat boundary metrics as in \cite[Lemma 2.2]{AF09}, so that $f$ is an isometry between the boundary components being identified. Apply the preceding flattening on each side. The two product collars glue smoothly to a compact flat product collar. All sectional curvatures outside these finitely many collars are strictly negative. Theorem~\ref{main}, with no domains $B_j$, applies.
\end{proof}

The same conclusion holds for the smooth manifold $W=D_{\widehat f}M$ constructed in \cite[Theorem 1.3]{AF09}. It also holds for the doubles in \cite[Theorem 1.6]{AF09}, by transporting an affine-double metric through the diffeomorphism supplied there. In Theorem 1.3 the conclusion concerns $W$ with its constructed smooth structure; the stated homeomorphism alone does not transport a smooth metric to a prescribed smooth structure on $D_fM$.

\section{The exotic modification}

\begin{corollary}\label{exotic}
The nonpositively curved metrics in the examples of \cite[Theorem 1.5]{AF09} can be chosen to have ergodic geodesic flow for Liouville measure.
\end{corollary}

\begin{proof}
Aravinda and Farrell first construct the manifold in the form
\[
 W\cong DM\mathbin{\#}\Sigma,
\]
where $\Sigma$ is a homotopy sphere, and then identify it diffeomorphically with a twisted double. The metric modification is supported in a compact region in the hyperbolic core, using the construction of Farrell and Jones \cite[Proposition 1.3]{FJ}. The cusp flattening is performed away from this region; see \cite[pp.~1056--1057]{AF09}.

Choose a slightly larger compact domain $B$ containing the entire metric modification, with boundary in the unchanged hyperbolic region. It is a disk with the local connected-sum modification in its interior, and is simply connected. It is disjoint from the cusp collars. Outside $B$ and those collars, the metric is strictly negatively curved, using the flattening of Section~2. Inside $B$ it is nonpositively curved by the cited construction. Theorem~\ref{main} therefore applies. Transporting the metric to the diffeomorphic twisted double preserves ergodicity.
\end{proof}

\section{Relation to earlier ergodicity results}

Wu proved that, on a closed surface of genus at least two with nonpositive curvature, finiteness of the connected components of $\{K<0\}$ implies that every flat geodesic is closed and that the geodesic flow is ergodic \cite[Theorem 1.6]{Wu15}. Here a flat geodesic is one along which the Gaussian curvature vanishes identically. Federico Rodriguez Hertz pointed out the connection with this work. Wu's proof excludes nonclosed flat geodesics by studying negatively curved components in an ideal triangle bounded by asymptotic geodesics. The argument here uses the prescribed product geometry of the collars to identify the singular set directly. In dimension two their torus sections are circles, and every singular geodesic is closed. In dimensions $d\geq3$, a flat torus section also contains nonclosed geodesics: on $\mathbb R^{d-1}/\Lambda$, the geodesic with velocity $v$ is closed precisely when $Tv\in\Lambda$ for some $T>0$. Such nonclosed singular geodesics still form part of the null set in \eqref{singular}.

Rodriguez Hertz, Rodriguez Hertz, and Wu \cite{RHW} prove ergodicity for closed surfaces of genus at least two with nonpositive curvature and a $C^{3+\mathrm{Lip}}$ metric, meaning that the third derivatives of the metric are locally Lipschitz. They also announce ergodicity for every closed connected smooth rank-one manifold of nonpositive sectional curvature in every dimension, with the proof to be presented separately.

Every metric considered in Theorem~\ref{main} has rank one: a geodesic passing through the nonempty region $U$, where all sectional curvatures are negative, admits no nonzero parallel normal Jacobi field. Consequently, the announced theorem implies the ergodicity conclusions of Theorem~\ref{main} and Corollaries~\ref{affine} and~\ref{exotic}. The present note supplies a direct proof for these metrics and the explicit description \eqref{singular} of their singular set. The Cartan--Hadamard argument also shows why the localized exotic modification introduces no additional singular geodesics.

Their surface proof establishes the measure-zero statement through Fourier decomposition, heat regularization, and commutator estimates, using the Guillemin--Kazhdan operators and Kannai's heat--wave formula. Here it follows from the fact that a geodesic confined to a flat product collar has zero normal velocity. The reduction to the singular set belongs to the classical regular-set ergodicity theory; a related geometric argument for flat neighborhoods appears in \cite[\S5]{CLMT}.

The examples in \cite{AF05} have strictly negative sectional curvature, although their underlying smooth manifolds admit no metric with nonpositive curvature operator. For their constructed negatively curved metrics, ergodicity follows directly from the classical theorem of Anosov \cite{Anosov}. Corollaries~\ref{affine} and~\ref{exotic} concern the particular choices of nonpositively curved metrics described above.

\section*{Acknowledgments}
I am grateful to Federico Rodriguez Hertz, Jana Rodriguez Hertz, and Weisheng Wu for the work that inspired this note. Their treatment of the singular set prompted the question of how the same strategy applies to the examples of Aravinda and Farrell. I warmly thank Jana and Federico for their encouraging comments, and Federico in particular for drawing my attention to Wu's earlier work \cite{Wu15}. I also acknowledge the geometric constructions of Aravinda and Farrell and the earlier work cited above, on which the argument depends.

I thank OpenAI for making ChatGPT available as a tool for mathematical discussion and writing. I used ChatGPT (OpenAI) to develop the geometric argument, examine its application to the Aravinda--Farrell constructions, check bibliographic information, and draft and revise the manuscript. This assistance included the description of singular geodesics in the product collars and the use of the Cartan--Hadamard theorem to exclude complete geodesics confined to the region of the exotic modification. I am responsible for the mathematical claims, exposition, and references.
\section*{Funding}
This work was supported by Proyecto PAPIIT, Direcci\'on General de Asuntos del Personal Acad\'emico, Universidad Nacional Aut\'onoma de M\'exico [grant number IN103324].

\begingroup
\small

\medskip
\noindent\textbf{Alberto Verjovsky}\\
Instituto de Matem\'aticas, Unidad Cuernavaca, Universidad Nacional Aut\'onoma de M\'exico.\\
Av. Universidad s/n, Col. Lomas de Chamilpa, C.P. 62210, Cuernavaca, Morelos, M\'exico.\\
\href{mailto:albertoverjovsky@gmail.com}{\nolinkurl{albertoverjovsky@gmail.com}}

\endgroup
\end{document}